%% file: manuscript.tex
\documentclass[opre,nonblindrev]{informs3}

\DoubleSpacedXI 

\usepackage{endnotes}
\let\footnote=\endnote
\let\enotesize=\normalsize
\def\notesname{Endnotes}%
\def\enoteformat{\rightskip0pt\leftskip0pt\parindent=1.275em
  \leavevmode\llap{\makeenmark}}

\usepackage{natbib}
 \bibpunct[, ]{(}{)}{,}{a}{}{,}%
 \def\bibfont{\small}%

\makeatletter
\def\@mb@citenamelist{cite,citep,citet,citealp,citealt,citepalias,citetalias}
\makeatother

\defcitealias{ch1formulation}{ARPA-E GO Comp. 2019}
\defcitealias{NERCStandards}{NERC 2020}

\TheoremsNumberedThrough     
\ECRepeatTheorems

\EquationsNumberedThrough    

\newcommand{\st}{\text{s.t.}}
\usepackage[ruled,vlined]{algorithm2e}
\SetKwFor{With}{with}{do}{end}
\usepackage[export]{adjustbox}

\begin{document}


\RUNAUTHOR{Petra and Aravena}

\RUNTITLE{Solving security-constrained OPF}

\TITLE{Solving realistic security-constrained optimal power flow problems}

\ARTICLEAUTHORS{%
\AUTHOR{Cosmin Petra}
\AFF{Lawrence Livermore National Laboratory, \EMAIL{petra1@llnl.gov}} 
\AUTHOR{Ignacio Aravena}
\AFF{Lawrence Livermore National Laboratory, \EMAIL{aravenasolis1@llnl.gov}}
} 

\ABSTRACT{%
%
We present a decomposition approach for obtaining good feasible solutions for the
security-constrained alternating-current optimal power flow (SCACOPF) problem at an industrial scale
and under real-world time and computational limits. The approach aims at complementing the existing
body of literature on bounding the problem via convex relaxations. It was designed for the participation in ARPA-E's Grid Optimization (GO) Competition Challenge 1. The challenge focused
on a near-real-time version of the SCACOPF problem where a base case operating point is optimized
taking into account possible single-element contingencies, after which the system adapts its
operating point following the response of automatic frequency drop controllers and voltage
regulators. Our solution approach for this problem relies on state-of-the-art nonlinear programming
algorithms and employs nonconvex relaxations for complementarity constraints, a specialized
two-stage decomposition technique with sparse approximations of recourse terms, and contingency
ranking and pre-screening.  The paper also outlines the
salient features of our implementation, such as fast model functions and derivatives evaluation, warm-starting
strategies, and asynchronous parallelism.  We discuss the results of the independent benchmark of our
approach done by ARPA-E's GO team in Challenge 1, which found that our methodology consistently produces high quality
solutions across a wide range of network sizes and difficulty. Finally, we conclude by outlining potential 
extensions and improvements of our methodology.
%
}%



\maketitle

%


%
%

\section{Motivation} \label{sec:motivation}

Industry practices for power systems scheduling and operation rely upon optimization tools that were
designed and specialized for conventional networks. In these networks, roughly speaking,
large-scale generators inject power into the bulk transmission system, then the transmission system
transport the power to step down substations, from which distribution grids draw power in order to
deliver it to final consumers. Emerging technologies disrupt this conventional model at several
levels, among them: (i) intermittent resources (wind and solar generation) make current security
and reliability requirements inadequate, (ii) distributed generation can reverse the direction of
power flows, and (iii) storage and demand response allow reshaping demand, moving it both in time
(batteries) and space (electric vehicles). Existing optimization tools do not acknowledge the
characteristics of these emerging technologies in modern power grids, leading to multiple problems
ranging from computational performance (solution time) to incomplete pricing, with potential
consequences that can affect from short-term security and resilience, up to long-term adequacy.

These challenges have inspired an important stream of research over the last decade, most
prominently in the study of convex relaxations for the optimal power flow problem (see
\citep{Molzahn2019} and references therein). Convex relaxations have nice theoretical properties. 
They offer a global bound on the solution quality and, whenever the solution to a convex relaxation
is found to be feasible for the original problem, such solution is guaranteed to be a global
optimum. At the same time, convex relaxations have a major practical disadvantage: they often have
solutions that are infeasible for the original problem, and reconstructing a feasible solution to
the original problem starting from the relaxation solution may be as complex as finding a good solution
from scratch. This practical disadvantage has prevented the direct use of convex relaxations by
power system operators.
Furthermore, such
feasible operating point should be \emph{secure}, in the sense that the system should be able to
continue operating under any \emph{credible contingency} event, without losing balance or falling
towards a blackout \citepalias{NERCStandards}. The optimization problem collecting all the
aforementioned considerations is known as the \emph{security-constrained alternating current optimal
power flow} (SCACOPF) problem and its importance cannot be overstated; it is solved (in the sense of
finding a feasible point) by almost every power system operator worldwide, on a daily basis, as part
of tasks ranging from short-term planning (unit commitment, economic dispatch) up to long-term
adequacy assessment (reliability studies, expansion planning). Solving SCACOPF accurately is
 a critical function upon which depend the reliability, security and efficiency of power
systems, as well as the correct functioning of other critical infrastructure dependent on
electricity.

In this paper we present a scalable method for computing good feasible solutions to the SCACOPF
problem. Our aim is to complement the recent literature on convex relaxations for the optimal power
flow problem, as well as, to present a solution approach that can be readily used by power system
operators in their daily operations.

The proposed method exploits state-of-the-art computational optimization and parallel computing,
while observing solution time requirements and computational resources currently available to system
operators. The method was designed within the scope of ARPA-E's Grid Optimization (GO) Competition
-- Challenge 1, part of a wider initiative by the U.S.  government agency to spur development of
new, modern and innovative software for the modernized grid. Our method was independently tested
during Challenge 1 using synthetic and real power grids, with sizes ranging from 1\,000 up to
30\,000 buses, that is, up to the size of the largest networks for which SCACOPF needs to be solved
worldwide. It provided satisfactory results in all cases, finding the best-known solution in most of
them and coming very close to the best-known solution in all other cases.

The rest of this paper is organized as follows. Section \ref{sec:state-of-the-art} reviews the
state-of-the-art in methods for SCACOPF with focus on methods capable of finding feasible
solution in a realistic setting. Section \ref{sec:scacopf_formulations} presents the SCACOPF
formulation used in Challenge 1 as well as our reformulation, approximations, and relaxations of
the complicating aspects of the problem. Section \ref{sec:decomposition} presents our solution approach,
based on decomposition ideas, relaxation and feasibility recovery, and sparse approximations for
recourse functions. Section \ref{sec:computational_implementation} outlines our computational
implementation and reviews important considerations when solving these problems in practice. Section
\ref{sec:numerical_results} summarizes our numerical experiments conducted during Challenge 1. Conclusions
 and directions for further research are presented in Section~\ref{sec:conclusions}.
 


\section{State-of-the-art and contributions} \label{sec:state-of-the-art}


The study of the optimal power flow problem (OPF) dates back, at least, to the 1960's, with the
first formulation of the problem being proposed in \citep{Carpentier1962}, following fear of
branch overloadings overlooked by simpler \emph{economic dispatch} procedures based on marginal generation
costs and loss coefficients \citep{Carpentier1979a}. The importance of security was recognized as
early as methods for this problem began to be developed; as J. Carpentier (France) noted ``\ldots
\emph{the general failure of 19 December 1973 lasted only 3h but caused a loss of production for the
country estimated to be at least the equivalent of 50 years' savings through economic dispatch}
\ldots'' (sic). This stark comparison between potential societal loses and energy cost savings,
especially when considering that the probability of failure of a typical power line over a year is
0.6\% \citep{Ekisheva2015}, motivated the study and the incorporation of security constraints in OPF
problems by the late 1970's.

The first approaches to SCACOPF were based on finding solutions to
the Karush-Kuhn-Tucker conditions of SCACOPF using diverse approximations for gradients and
Hessians within first- and second-order methods \citep{Cory1972, Alsac1974}. These approaches proven
effective for small systems, but were unable to scale to utility-scale power grids. Other methods
based on simplification of the power flow equations, such as their linear \citep{Stott1978} or
quadratic approximation \citep{Carpentier1979b}, and sequential evaluation of base case and
contingency subproblems, achieve better scalability at the cost of accuracy and suboptimality.
Nevertheless, the comparative advantages of the latter methods with respect to existing industry
tools for near-real-time generator dispatch (at the time) led to their early adoption by a few power
system operators, such as the New York Power Pool \citep{Elacqua1982}.

Subsequent advances in the area during the 80's came from the adoption of formal
decomposition techniques in SCACOPF, which permitted more detailed SCACOPF models. In particular,
Benders decomposition was employed to solve a simplified SCACOPF (linearized, real power only) with
corrective actions (\textit{i.e.},  dispatch decisions in contingency subproblems) exactly
\citep{Pereira1985}, a work that was later extended to the full non-linear SCACOPF problem using the
generalized Benders decomposition \citep{Monticelli1987}, though this extension does not guarantee
convergence because of non-convexity of contingency subproblems. During this decade most industrial
implementations were based on linear programming technology \citep{Alsac1990}, which was much more
mature than nonlinear programming (NLP) technology at that point. Only by the end of the decade,
non-linear programming techniques became competitive alternatives to construct SCACOPF algorithms that
could scale up to the requirements of the industry \citep{Papalexopoulos1989}.

Despite continuing progress in the subsequent years, by the mid-90's there were still several
challenges that needed to be addressed both in the formulation and solution of the SCACOPF problem
\citep{Momoh1997}. One of these challenges is the nonconvex nature of the problem, which prevents its 
direct use in computing electricity prices (dual variables of the problem) and prevents guaranteeing
global optimality. Operators have resorted to convex (mostly linear) approximations of SCACOPF to
compute prices and to obtain global solutions for the approximation; recovering feasible solutions
for the original SCACOPF problem in postprocessing steps. 
The nonconvexity challenge, however, remains unsolved and an area of intensive research until now.
Other challenges observed were mainly computational: incorporating additional restrictions into the
the problem, scaling to larger interconnected networks, and improving solution speed for the use of
SCACOPF engines in near-real-time operations. These additional restrictions for SCACOPF correspond
to either (\textit{i}) more detailed models for the physics of the power grid, such as transient
stability requirements for the transition between base and contingencies, as later captured in
\citep{Bruno2002} and \citep{Yuan2003}, or (\textit{ii}) restrictions imposed by the limitations of
control systems and system operators, such as a limit on the number of corrective actions that can
be taken between the base case and any given contingency, as studied later in \citep{Capitanescu2011}
and \citep{Phan2015}.

Parallel computing has demonstrated promising results for solving the aforementioned scalability and
speed challenges of SCACOPF in a systematic fashion. The first parallel computing approach to
SCACOPF \citep{Rodrigues1994} used a linear approximation of the problem and it was based on
decomposing the computations required to run a dual simplex method with constraint generation.
In~\cite{Qiu2005} and ~\cite{petra_14_augIncomplete,petra_14_realtime} the nonlinear SCACOPF problem is solved in parallel by decomposing the
linear algebra of interior point methods. More recent parallel computing
approaches, such as \citep{Liu2013}, decompose the SCACOPF problem at the level of the formulation,
usually separating the problem into base case and contingency subproblems. Subproblems are solved
using \emph{off-the-shelf} nonlinear programming packages such as Ipopt \citep{Watcher2006}, and
coordination among subproblems is achieved through a first-order optimization method, in a very
similar fashion to earlier work in the area \citep{Pereira1985}.

Two key observations can further improve the scalability and speed of virtually any SCACOPF engine.
First, in real instances, only a handful of the (possibly) thousands of contingencies of a system
are severe enough to induce changes to the base case solution. Typical SCACOPF industrial
implementations will rank contingencies based on certain \emph{severity index} and will optimize
only considering the contingencies above a certain threshold, an approach that depends on
well-defined indexes and requires significant parameter tuning. A different idea is to search 
directly for a subset of contingencies such that solving SCACOPF with the subset is (almost)
equivalent to solving the full SCACOPF problem. \cite{Bouffard2005} proposes to find such subset by
looking at dual values for balance constraints in contingencies, whereas \citep{Capitanescu2007}
proposes to do so looking at the primal violations incurred in contingency constraints when ignoring
them. The second observation is that, in large real networks, contingencies mostly affect nearby
elements, while most of the network is not severely disrupted. \cite{Karoui2008} exploits this
observation to develop approximation models for contingency subproblems, where areas nearby the
failed element(s) are modeled in full detail, whereas the rest of the system is reduced to radial
equivalent circuits, effectively reducing the size and complexity of contingency subproblems.
Parallel computing, contingency subset selection and contingency subproblem approximation were later
use in tandem in the algorithm proposed in \citep{Platbrood2014}, which was able to find a good
feasible solution to the problem for a 9\,241-bus system with 12\,000 contingencies in 1 hour.

In parallel to these developments, second-order cone programming \citep{Jabr2006} and semidefinite
programming \citep{Bai2008} relaxations for the optimal power flow problem started to be
investigated as potential alternatives to solve the inability of existing methods to achieve
globally optimal solutions and compute equilibrium prices (whenever they exist). The approach show
promise, with many of the usual (small) test cases of the literature being solved globally for these
technique for the first time in \citep{Lavaei2011}. This work also presented sufficient conditions
for achieving zero duality gap in general networks.  Several new relaxations and strengthening
techniques for existing relaxations have been proposed over the years, and the area remains a very
active field of research \citep{Molzahn2019}. Convex relaxations, however, have rarely been studied
in the context of full SCACOPF problems. In one of theses studies, \cite{Madani2016} evaluates the
use of the semidefinite programming relaxation for solving SCACOPF and proposes a penalization
technique to recover feasible solutions whenever the rank-1 condition is not satisfied. A major
shortcoming of the approach, as of any known convex relaxation of OPF, lies on the necessity of this
feasibility recovery step, which loses the global optimality guarantees of the relaxation, and that
is not trivial to carry out, in general, without moving far away from the solution found by the
relaxation. These concerns have so far prevented the deployment of convex relaxations of SCACOPF in
industrial implementations. 

Recent studies in SCACOPF have focused on improving the representation and tractability of realistic
coupling constraints between the base case and contingencies \citep{Dvorkin2018, Velloso2020}, which
permits not only optimizing the base case solution to make it secure, but also optimize the primary
response of generators to contingencies. At the same time, novel machine learning techniques have
begun to be proposed and theoretically validated for SCACOPF \citep{Venzke2021}, opening the door
for future progress in the area, as is required by the evolving needs of power grids.


Our contributions to this vast body of literature are twofold. We first propose
an efficient decomposition approach to solve SCACOPF in near-real-time for real-sized power systems,
considering the action of frequency drop control and voltage regulators as recourse actions
\citepalias{ch1formulation}. Our approach employs relaxations and feasibility recovery in order to
handle the non-smoothness of automatic controller actions, and a sparse and an optimistic 
approximation for recourse terms corresponding to penalizations for infeasibilities in contingency
subproblems. It differs from previous approaches in that it does not require tens of passes over the
list of contingencies, such as \citep{Pereira1985}, \cite{Liu2013}, but only a 2 or 3, which is as
much as an operator can afford during near-real-time operations. 

Second, we detail an implementation of our approach using state-of-the-art nonlinear programming and
parallel computing techniques, which permit the solution of full ACOPF problems for all
contingencies, as opposed to the compressions used in the literature \citep{Karoui2008,
Platbrood2014}, which can lead to misrepresentation of the effects of contingencies of the system
(namely, global effects of contingencies due to automatic action of drop frequency controllers). Our
methodology allows us to effectively compute good SCACOPF solutions for systems with up to 70\,000
buses and 20\,0000 contingencies in 15 minutes, that is, its capabilities exceed the needs of any
current system operator worldwide, solving the computational challenges required to bring SCACOPF to
real time operations.

%

\section{SCACOPF formulations} \label{sec:scacopf_formulations}

Challenge 1 posed a near-real-time version of SCACOPF \citepalias{ch1formulation}: finding the
lowest cost generator production levels that would satisfy consumer demand for a base case (where
all power equipment is available), so that operational constraints can be satisfied even if a
contingency (sudden loss of certain equipment) were to happen, while taking into account the
actions of automatic controllers, \textit{e.g.}, frequency drop control and voltage regulators, in response to
contingencies. In this section we first present an stylized version of such formulation of
SCACOPF\endnote{Specifically, we do not present the models for transformer branches and switched
shunts, the models for piecewise linear costs and penalty functions, and we do not consider multiple
areas when presenting the model for frequency drop control, nor the case where these controllers may
have negative constants. The complete model can be found in \citepalias[Sec. 3]{ch1formulation}.},
highlighting its most important features and complicating constraints. Then, we present and justify
our smooth reformulation of the problem, with particular attention to the modeling of automatic
controllers. Regarding notation, unless indicated otherwise, we use lowercase letters for variables,
uppercase letters for parameters, and calligraphic letters for sets.

\subsection{ARPA-E GO Challenge 1 official SCACOPF formulation} \label{sec:arpa-e_formulation}

Let $\mathcal{K} = \{k_0, k_1, \ldots, k_K\}$ be the set of indices of power flow cases in SCACOPF,
$k_0$ corresponding to the base case and $k_1, \ldots, k_K$ corresponding to contingencies.
We consider that each contingency correspond to either the loss of a single generating unit or a
single power line, a criterion known as N-1 security in technical literature. In this setting, we
formulate SACOPF as follows.

Constraints \eqref{eq:voltage_bounds} impose bounds on bus (node) voltages $v_{n,k}$ at all buses
$n \in \mathcal{N}$,
\begin{equation}
	\underline{V}_{n,k} \leq v_{n,k} \leq \overline{V}_{n,k} ~~ \forall n \in \mathcal{N}, k \in
\mathcal{K}, \label{eq:voltage_bounds}
\end{equation}
which can correspond to either normal operation bounds for the base case $k = k_0$, or emergency
bounds for contingencies $k \neq k_0$, with $\underline{V}_{n,k} \leq \underline{V}_{n, k_0}$ and
$\overline{V}_{n, k_0} \leq \overline{V}_{n, k}$ for any $n \in \mathcal{N}, k \in
\mathcal{K} \setminus \{k_0\}$.

Constraints \eqref{eq:p_bounds} and \eqref{eq:q_bounds} impose bounds on active power injection
$p_{g,k}$ and reactive power injection $q_{g,k}$ for all generators $g \in \mathcal{G}^k$,
\begin{gather}
	\underline{P}_g \leq p_{g,k} \leq \overline{P}_g ~~ \forall g \in \mathcal{G}^k, k \in
\mathcal{K}, \label{eq:p_bounds} \\
	\underline{Q}_g \leq q_{g,k} \leq \overline{Q}_g ~~ \forall g \in \mathcal{G}^k, k \in
\mathcal{K}, \label{eq:q_bounds}
\end{gather}
where $\mathcal{G}^k \subseteq \mathcal{G}^{k_0}$ for any $k \in \mathcal{K} \setminus \{k_0\}$,
that is, generators available in a contingency $k$ are the same as those available in the base case,
minus the failed generators.

Constraints \eqref{eq:p_branch} -- \eqref{eq:capacity_branch} model the power flows entering the
origin terminal $o$ of branch (edge) $e \in \mathcal{E}^k$ in terms its series conductance $G_e$, series
susceptance $B_e$, charging susceptance $B^{CH}_e$, thermal rating $R^k_e$, and origin $o(e) \in
\mathcal{N}$ and destination $d(e) \in \mathcal{N}$ bus voltage magnitudes $v_{o(e),k}, v_{d(e),k}$
and angles $\theta_{o(e),k}, \theta_{d(e),k}$, in power flow case $k \in \mathcal{K}$, in the
following fashion: 
\begin{align}
	p_{e,o,k} ~&= G_e v_{o(e),k} - G_e \, \cos\big(\theta_{o(e),k} - \theta_{d(e),k}\big) \,
                                          v_{o(e),k} \, v_{d(e),k} ~- \nonumber \\
		     & \quad B_e \, \sin\big(\theta_{o(e),k} - \theta_{d(e),k}\big) \, 
                            v_{o(e),k} \, v_{d(e),k} ~~ \forall e \in \mathcal{E}^k, k \in \mathcal{K},
                            \label{eq:p_branch} \\
	q_{e,o,k} ~&= -\big(B_e + B_e^{CH}/2\big) \, v_{o(e),k}^2 + 
                  B_e \, \cos\big(\theta_{o(e),k} - \theta_{d(e),k}\big) \,
                         v_{o(e),k} \, v_{d(e),k} ~- \nonumber \\
		     & \quad G_e \, \sin\big(\theta_{o(e),k} - \theta_{d(e),k}\big) \,
                            v_{o(e),k} \, v_{d(e),k} ~~ \forall e \in \mathcal{E}^k, k \in \mathcal{K},
                            \label{eq:q_branch}
\end{align}
\begin{equation}
	\sqrt{p_{e, o, k}^2 + q_{e, o, k}^2} \leq R_{e,k} v_{o(e), k} + \sigma_{e, o, k}, ~
		\sigma_{e, o, k} \geq 0 ~~ \forall e \in \mathcal{E}^k, k \in \mathcal{K}. \label{eq:capacity_branch}
\end{equation}
Constraints \eqref{eq:p_branch} define the active power flow $p_{e,o,k}$ entering terminal $o$,
constraints \eqref{eq:q_branch} define the reactive power flow $q_{e,o,k}$ entering terminal $o$,
and constraints \eqref{eq:capacity_branch} enforces the thermal limits of the branch, allowing for
violations through the slack variable $\sigma_{e, o, k}$. 
Similar to the case of generators, we have that $\mathcal{E}^k \subseteq \mathcal{E}^{k_0}$, and
similar to the case of buses, $R_{e,k} \geq R_{e,k_0}$, for any $k \in \mathcal{K} \setminus
\{k_0\}$. Power flows entering the destination terminal $d$ are defined in an analogous fashion,
\textit{i.e.}, interchanging $o$ and $d$ in constraints \eqref{eq:p_branch} --
\eqref{eq:capacity_branch}.

Constraints \eqref{eq:p_balance} -- \eqref{eq:balance_slacks_nn}, below, model the active and reactive power
balance at each bus of the grid $n \in \mathcal{N}$, under every power flow case $k \in
\mathcal{K}$.
\begin{gather}
	\sum_{g \in \mathcal{G}^k(n)} p_{g,k} - P^D_n - G^{SH}_n v_{n,k}^2 -
		\sum_{\substack{e \in \mathcal{E}^k: \\ o(e) = n}} p_{e,o,k} - 
		\sum_{\substack{e \in \mathcal{E}^k: \\ d(e) = n}} p_{e,d,k} =
        \sigma^{P+}_{n, k} - \sigma^{P-}_{n, k} ~~ \forall n \in \mathcal{N}, k \in \mathcal{K}
        \label{eq:p_balance} \\
	\sum_{g \in \mathcal{G}^k(n)} q_{g,k} - Q^D_n + B^{SH}_n v_{n,k}^2 -
		\sum_{\substack{e \in \mathcal{E}^k: \\ o(e) = n}} q_{e,o,k} - 
		\sum_{\substack{e \in \mathcal{E}^k: \\ d(e) = n}} q_{e,d,k} = 
        \sigma^{Q+}_{n, k} - \sigma^{Q-}_{n, k} ~~ \forall n \in \mathcal{N}, k \in \mathcal{K} 
        \label{eq:q_balance}\\
	\sigma^{P+}_{n, k}, \sigma^{P-}_{n, k}, \sigma^{Q+}_{n, k}, \sigma^{Q-}_{n, k} \geq 0 ~~ 
        \forall n \in N, k \in \mathcal{K} \label{eq:balance_slacks_nn}
\end{gather}
Here, active power balance constraints \eqref{eq:p_balance} account for (\textit{i}) active power
injections of all generators $g$ at bus $n$, denoted as $\mathcal{G}^k(n)$, (\textit{ii}) active
demand from consumers $P^D_n$, (\textit{iii}) shunt conductance $G^{SH}_n$, and (\textit{iv}) power
injected to branches $e \in \mathcal{E}^k$ whose origin or destination terminals connect to bus $n$.
Violations of power balance are allowed through the slacks $\sigma^{P+}_{n,k}$ and
$\sigma^{P-}_{n,k}$, which are forced to be nonnegative by constraints \eqref{eq:balance_slacks_nn}.
Reactive power balance constraints \eqref{eq:q_balance} are defined in an analogous fashion, with
$Q^D_n$ corresponding to reactive demand from consumers and $B^{SH}_n$ corresponding to shunt
admittance, at bus $n$.

In the near-real-time version of SCACOPF of \cite{ch1formulation}, contingencies originate on
external events to the system (\textit{e.g.}, a tree falling down on a power line) while the system is
operated in its base case condition. Following the contingency event, the power system adjusts --
departing from its base case condition -- according to the actions taken by automatic control
systems (frequency drop control and voltage regulators) until reaching a new quasi-stationary point
which we model as the contingency power flow constraints presented above (constraints
\eqref{eq:voltage_bounds} -- \eqref{eq:balance_slacks_nn} for $k \in \mathcal{K} \setminus
\{k_0\}$). This dynamic adjustment process couples the base case conditions with the contingency
conditions. The SCACOPF model captures this coupling as follows.

The frequency drop control constraints following a contingency event are  modeled as:
\begin{gather}
	p_{g,k} + \rho^+_{g,k} - \rho^-_{g,k} = p_{g,k_0} + A_{g} \delta_{k} ~~
        \forall g \in \mathcal{G}^k, k \in \mathcal{K} \setminus \{k_0\} ,\label{eq:drop_control}\\
	0 \leq \rho^+_{g,k} \perp \overline{P}_g - p_{g,k} \geq 0 ~~ 
        \forall g \in \mathcal{G}^k, k \in \mathcal{K} \setminus \{k_0\},
        \label{eq:drop_control_up_saturation} \\
	0 \leq \rho^-_{g,k} \perp p_{g,k} - \underline{P}_g \geq 0 ~~
        \forall g \in \mathcal{G}^k, k \in \mathcal{K} \setminus \{k_0\}.
        \label{eq:drop_control_down_saturation}
\end{gather}
Given a contingency $k \in \mathcal{K} \setminus \{k_0\}$, the active power injected by generator $g
\in \mathcal{G}^k$ to the power system is proportional to a system-wide signal $\delta_k$ (emulating
frequency deviation), with proportionality constant $A_g \geq 0$. This proportional control action is
saturated at the technical limits for production of generator $g$, which is accomplished by the use
of the slack variables $\rho^+_{g,k}$ and $\rho^-_{g,k}$ (corresponding to over saturation above
and below), and the complementarity constraints \eqref{eq:drop_control_up_saturation} and
\eqref{eq:drop_control_down_saturation}, which model saturation at the upper and lower limits,
respectively.

Similarly, the behavior of regulators of generators following a contingency event is
modeled as:
\begin{gather}
	\nu^+_{n,k} - \nu^-_{n,k} = v_{n,k} - v_{n, k_0} ~~ 
        \forall n \in \mathcal{N}^k, k \in \mathcal{K} \setminus \{k_0\} 
        \label{eq:v_regulator_diff} \\
	0 \leq \nu^-_{n,k} \perp \overline{Q}_g - q_{g,k} \geq 0 ~~ 
        \forall g \in \mathcal{G}^k(n), n \in \mathcal{N}^k, k \in \mathcal{K} \setminus \{k_0\} 
        \label{eq:v_regulator_up} \\
	0 \leq \nu^+_{n,k} \perp q_{g,k} - \underline{Q}_{g} \geq 0 ~~ 
        \forall g \in \mathcal{G}^k(n), n \in \mathcal{N}^k, k \in \mathcal{K} \setminus \{k_0\},
        \label{eq:v_regulator_down}
\end{gather}
where $\mathcal{N}^k := \{n \in \mathcal{N} \mid \mathcal{G}^k(n) \neq \emptyset\}$ is the set of
buses with generators connected to them (where voltage is controlled).  In a contingency event $k
\in \mathcal{K} \setminus \{k_0\}$, generator $g \in \mathcal{G}^k(n)$ regulates its active power
injection trying to maintain the base case voltage magnitude at bus $n$. If the voltage cannot be
maintained, and there is decrease (increase) in the contingency voltage, the generator must be
injecting its maximum (minimum) reactive power to the network. This saturated control actions is
modeled using the slacks variables $\nu^-_{n,k}$ and $\nu^+_{n,k}$ (corresponding to voltage
decreases and increases in contingency $k$ with respect to the base case), and the complementarity
constraints \eqref{eq:v_regulator_up} and \eqref{eq:v_regulator_down}, which model reactive
injection in cases of voltage decrease and increase, respectively.

The objective function of SCACOPF can be cast as equation \eqref{eq:objective},
\begin{equation}
    \begin{aligned}
        \min_{\substack{p,q,v,\theta \\ \nu, \rho}} ~~
            & \sum_{g \in \mathcal{G}^{k_0}} f_g(p_{g,k_0}) + 
                \sum_{e \in \mathcal{E}^{k_0}} \left(\gamma^S(\sigma_{e,o,k_0}) + 
                                                     \gamma^S(\sigma_{e,o,k_0})\right) + \\
            & \sum_{n \in \mathcal{N}} 
                \left(\gamma^P(\sigma^{P+}_{n,k_0}) + \gamma^P(\sigma^{P-}_{n,k_0}) +
                      \gamma^Q(\sigma^{Q+}_{n,k_0}) + \gamma^Q(\sigma^{Q-}_{n,k_0})\right) + \\
            & \frac{1}{|\mathcal{K}| - 1} ~ \sum_{k \in \mathcal{K} \setminus \{k_0\}} \Bigg(
                \sum_{e \in \mathcal{E}^k} \left(\gamma^S(\sigma_{e,o,k}) + 
                                                 \gamma^S(\sigma_{e,o,k})\right) + \\
            & \qquad \sum_{n \in \mathcal{N}} 
                        \left(\gamma^P(\sigma^{P+}_{n,k}) + \gamma^P(\sigma^{P-}_{n,k}) +
                              \gamma^Q(\sigma^{Q+}_{n,k}) + \gamma^Q(\sigma^{Q-}_{n,k})\right)
              \Bigg)
    \end{aligned}
    \label{eq:objective}
\end{equation}
where $f_g$ is the piecewise linear convex cost function of generator $g \in \mathcal{G}^{k_0}$,
$\gamma^S$ is the piecewise linear convex penalty function that determines how much it costs to
overload a branch, and $\gamma^P$ and $\gamma^Q$ are the piecewise linear convex penalty functions
that determine how much it costs to incur in active ($P$) and reactive ($Q$) power imbalances. The
objective function, thus, accounts for generators' cost, overload and imbalance penalties in the
base case, and overload and imbalance penalties in contingencies.

\subsection{SCACOPF (approximate) re-formulation} \label{sec:gollnlp_formulation}

Our approach to SCACOPF is based on decomposing the problem \eqref{eq:voltage_bounds} --
\eqref{eq:objective} into base case and contingency subproblems, iterating between them until an
acceptable solution is found. As such, our approach requires efficient methods for solving both the
base case and contingency subproblems, for which we resort to interior point methods
\citep{Watcher2006}. This design choice imposes the limitation that we cannot include (directly)
non-smooth expressions in our formulation. The model in the previous section includes non-smooth
expressions in constraints \eqref{eq:capacity_branch}, \eqref{eq:drop_control_up_saturation},
\eqref{eq:drop_control_down_saturation}, \eqref{eq:v_regulator_up} and \eqref{eq:v_regulator_down}.
At the same time, in order to meet the performance necessary for our method to be effective, we need
to simplify every aspect of the formulation that can be simplified without compromising solution
quality. In what follows we detail and justify the reformulations, relaxations and approximations to
problem \eqref{eq:voltage_bounds} -- \eqref{eq:objective} used in our approach to deal with
the aforementioned challenges.


We reformulate the non-smooth thermal limits constraints \eqref{eq:capacity_branch} simply squaring
both sides of the inequality, which yields
\begin{equation}
	p_{e, o, k}^2 + q_{e, o, k}^2 \leq (R_{e,k} v_{o(e), k} + \sigma_{e, o, k})^2 ~~ 
        \forall e \in \mathcal{E}^k, k \in \mathcal{K}.
        \label{eq:capacity_branch_squared}
\end{equation}
Whereas constraint \eqref{eq:capacity_branch_squared} is nonconvex, it should be noted that when
used within a log-barrier algorithm, it will be added to the objective of each barrier
subproblem
as
\begin{equation*}
    \log\left((R_{e,k} v_{o(e), k} + \sigma_{e, o, k})^2 - p_{e, o, k}^2 - q_{e, o, k}^2\right)
\end{equation*}
which is a convex expression and, in fact, a self-concordant barrier function for the second-order
cone \citep{Boyd2004}.

\begin{figure}
\centering
\begin{minipage}{.5\textwidth}
    \centering
    {\small (a)}
\end{minipage}%
\begin{minipage}{0.5\textwidth}
    \centering
    {\small (b)}
\end{minipage}
\begin{minipage}[t]{.5\textwidth}
    \centering
    \includegraphics[scale=0.75]{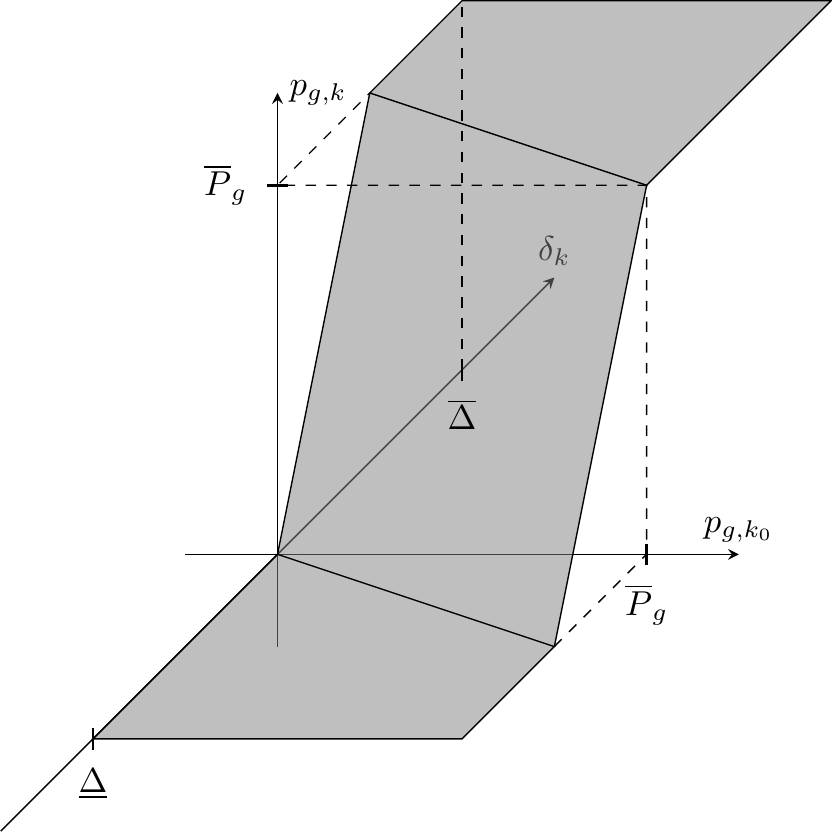}
\end{minipage}%
\begin{minipage}[t]{0.5\textwidth}
    \centering
    \includegraphics[scale=0.75]{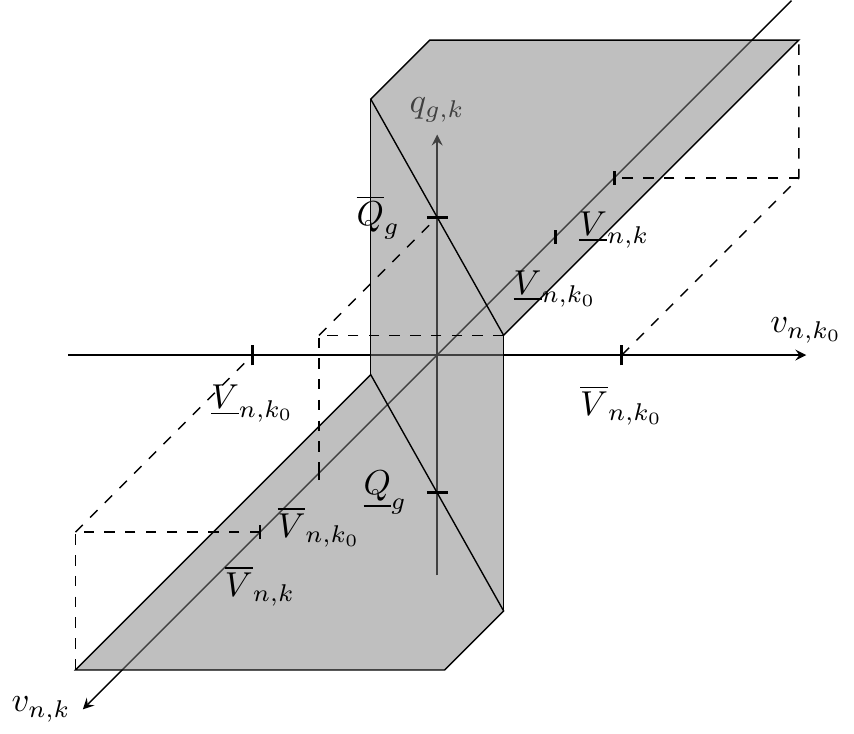}
\end{minipage}
\caption{Feasible sets of coupling constraints. (a) corresponds to the feasible set
$\mathcal{C}^p_{g,k}$ of constraints \eqref{eq:drop_control} --
\eqref{eq:drop_control_down_saturation} in the space of $p_{g,k_0}$, $\delta_k$, and $p_{g,k}$. (b)
corresponds to the feasible set $\mathcal{C}^v_{n,g,k}$ of constraints \eqref{eq:v_regulator_diff}
-- \eqref{eq:v_regulator_down} in the space of $v_{n,k_0}$, $v_{n,k}$, and $q_{g,k}$, plotted with
the origin at $(v_{n,k_0}=1, v_{n,k}=1, q_{g,k}=0)$ for clarity.}
\label{fig:feasible_sets_coupling}
\end{figure}

The non-smoothness due to complementarity relations in constraints \eqref{eq:drop_control_up_saturation},
\eqref{eq:drop_control_down_saturation}, \eqref{eq:v_regulator_up}, and \eqref{eq:v_regulator_down}
pose a  more considerable challenge. Existing smooth formulations, such as the Fischer–Burmeister
function \citep{Chen2008} -- which replaces a complementarity relation $0 \leq x \perp y \geq 0$ by
$\phi(x, y) = \sqrt{x^2 + y^2} - x - y = 0$, or by its squared (smooth) version $\phi(x, y)^2 = 0$ -- were
initially probed as an alternative but were ultimately abandoned due to frequent numerical issues
when optimizing over near $(0, 0)$. Convex relaxations of these constraints were also
investigated. In particular, let us define the sets
\begin{align*}
    \mathcal{C}^p_{g,k} &= \{(p_{g,k_0}, p_{g,k}, \delta_k) \in 
                             [\underline{P}_g, \overline{P}_g] \times 
                             [\underline{P}_g, \overline{P}_g] \times
                             [\underline{\Delta}, \overline{\Delta}] \mid 
                             \exists (\rho^+_{g,k}, \rho^-_{g,k}) \in \mathbb{R}^2_+: 
                             \eqref{eq:drop_control}, \eqref{eq:drop_control_up_saturation}, 
                             \eqref{eq:drop_control_down_saturation}\}, \\
    \mathcal{C}^v_{n,g,k} &= \{(v_{n,k_0}, v_{n,k}, q_{g,k}) \in 
                             [\underline{V}_{n,k_0}, \overline{V}_{n,k_0}] \times 
                             [\underline{V}_{n,k}, \overline{V}_{n,k}] \times
                             [\underline{Q}_g, \overline{Q}_g] \mid 
                             \exists (\nu^+_{n,k}, \nu^-_{n,k}) \in \mathbb{R}^2_+: 
                             \eqref{eq:v_regulator_diff}, \eqref{eq:v_regulator_up}, 
                             \eqref{eq:v_regulator_down}\},
\end{align*}
where $\underline{\Delta}$ and $\overline{\Delta}$ are bounds on $\delta_k$ for any $k$ (which can
be computed via interval algebra, for example). Sketches for these sets are presented in Fig.
\ref{fig:feasible_sets_coupling}. We note that both $\mathcal{C}^p_{g,k}$ and
$\mathcal{C}^v_{n,g,k}$ have polyhedral convex hulls, which can be easily described as convex
combinations of, at most, 8 extreme points. This approach, while numerically stable, leads to very
loose relaxations, for which there no incentive in the objective for the optimal to remain near a
feasible point and was ultimately abandoned due to the impossibility of recovering good feasible
solutions from it. 

Finally, we opted for relaxing complementarity relations in constraints
\eqref{eq:drop_control_up_saturation}, \eqref{eq:drop_control_down_saturation},
\eqref{eq:v_regulator_up} and \eqref{eq:v_regulator_down} by replacing $0 \leq x \perp y \geq 0$ by
    $x \geq 0$, $y \geq 0$, and $x \cdot y \leq \epsilon \cdot (\overline{X} - \underline{X}) \cdot 
                                                      (\overline{Y} - \underline{Y})$,
where $\epsilon$ is a small constant (\textit{e.g.},  $\epsilon = 10^{-4}$) and $\underline{X},
\overline{X}$ and $\underline{Y}, \overline{Y}$ are bounds for $x$ and $y$. This approach can be
thought as the middle ground between ill-conditioned smooth equivalents and well-conditioned but
loose convex approximations. The non-convex relaxation keeps numerical issues under control (in the
worst case, one can simply increase $\epsilon$), while providing almost-feasible solutions, suitable
for approximate evaluations within our iterative decomposition procedure for SCACOPF.

Regarding simplifications applied to the formulation of SCACOPF, we note that piecewise linear
penalty functions $\gamma^S, \gamma^P$, and $\gamma^Q$, nonetheless simple from a theoretical
perspective, have a detrimental performance impact, which is due to two main factors. First, these
terms are abundant (they are required multiple times for every bus and branch in the power grid) and
each of them requires an auxiliary variable and handful of linear constraints to be expressed in a
smooth fashion, significantly enlarging the size of the problem to what it would be its size
otherwise (approximately by 30\%). Secondly, the large slopes of these penalty functions towards
large slacks (5 orders of magnitude larger than other terms in the gradient of the objective
function and the Jacobian of the constraints \citep{ch1formulation}) can cause numerical issues
whenever they are active (have a positive value). 
Following these observations, we replaced these penalties functions by quadratic approximations of
the form $\tilde{\gamma}(x) = a_1 x + a_2 x^2$, where $a_1$ and $a_2$ are computed so as to match
the slope of the $\gamma$ functions at the origin and at the beginning of their second bin
(corresponding to small slacks). These approximations can be placed directly in the objective avoiding
extra variables and constraints. At the same time, our selection of their coefficients prevents the
aforementioned numerical issues, without degradation of solution quality in our computational 
experiments.


\section{Decomposition approach} \label{sec:decomposition}

We design our decomposition approach to SCACOPF with the computational and time limitations present
in real power grid operations, which were emulated in the GO Competition -- Challenge 1: up to 144
parallel processors and either 10 minutes (near-real-time operations) or 45 minutes (offline)
maximum elapsed time. Under these restrictions, we can only afford to evaluate contingencies (in
parallel) a couple of times (less than 10 for most realistic systems) in order to derive preventive
actions in the base case to ensure security constraints are met. These limiting factors on the
number of iterations we can perform, along with the non-convexity of subproblems, drove us away from
traditional decomposition approaches such as (generalized) Benders and ADMM methods, which typically
require tens of iterations to arrive at a good solution. Instead, we develop a specialized approach
that derives preventive actions to meet security constraints based on computationally affordable approximations for contingency
penalty terms (third and fourth lines in the objective function \eqref{eq:objective}). We
emphasize that our approach is designed to provide good feasible solutions, which do not necessarily 
satisfy accurately the (local) optimality conditions for problem \eqref{eq:voltage_bounds} --
\eqref{eq:objective}. Solutions produced by our approach can later be compared against convex
relaxations of \eqref{eq:voltage_bounds} -- \eqref{eq:objective} (\textit{e.g.}, second-order cone relaxation
of base case only) if a quality guarantee is required.

This section presents our decomposition approach for SCACOPF in terms of the subproblems involved,
feasibility recovery for contingencies, update rules and our contingency pre-screening methodology,
which will become the building blocks for our SCACOPF solver. The particular assembly of these
blocks in the asynchronous parallel computing implementation of our SCACOPF solver is discussed
later in section \ref{sec:computational_implementation}.


\subsection{Subproblems' formulation} \label{sec:subproblems}

We decompose \eqref{eq:voltage_bounds} -- \eqref{eq:objective} -- with the approximations described
in the preceding section -- by stages, leading to a first-stage master subproblem, corresponding to
the base case plus penalties coming from the contingencies, and many second-stage subproblems,
corresponding to contingencies. We formulate the \textit{first-stage master subproblem} as
\begin{equation}
\begin{aligned}
        \min_{p,q,v,\theta,\sigma} ~~
            & \sum_{g \in \mathcal{G}^{k_0}} f_g(p_{g,k_0}) + 
                \sum_{e \in \mathcal{E}^{k_0}} \left(\tilde{\gamma}^S(\sigma_{e,o,k_0}) + 
                                                     \tilde{\gamma}^S(\sigma_{e,o,k_0})\right) + \\
            & \sum_{n \in \mathcal{N}} 
                \left(\tilde{\gamma}^P(\sigma^{P+}_{n,k_0}) + \tilde{\gamma}^P(\sigma^{P-}_{n,k_0}) +
                      \tilde{\gamma}^Q(\sigma^{Q+}_{n,k_0}) + \tilde{\gamma}^Q(\sigma^{Q-}_{n,k_0})\right) + \\
            & \frac{1}{|\mathcal{K}| - 1} ~ \sum_{k \in \mathcal{K} \setminus \{k_0\}}
                r_k(\boldsymbol{p}_{k_0}, \boldsymbol{v}_{k_0}) \\
        \st \quad&  \eqref{eq:voltage_bounds}\,\text{--}\,\eqref{eq:balance_slacks_nn}: k = k_0.
\end{aligned} \label{eq:master}
\end{equation}
where $r_k$ are recourse functions corresponding to the objective contributions of second-stage
subproblems, \textit{i.e.}, penalties for not abiding with security constraints. Furthermore, the
second-stage contingency subproblems that define the functions $r_k$ are formulated as 
\begin{equation}
    \begin{aligned}
        r_j(\boldsymbol{p}_{k_0}, \boldsymbol{v}_{k_0}) = \quad & \\
            \min_{\substack{p,q,v,\theta \\ \sigma, \nu, \rho}} ~~& 
            \sum_{e \in \mathcal{E}^j}
                \left(\tilde{\gamma}^S(\sigma_{e,o,j}) + \tilde{\gamma}^S(\sigma_{e,o,j})\right) + \\
            & \sum_{n \in \mathcal{N}} 
                \left(\tilde{\gamma}^P(\sigma^{P+}_{n,j}) + \tilde{\gamma}^P(\sigma^{P-}_{n,j}) +
                      \tilde{\gamma}^Q(\sigma^{Q+}_{n,j}) + \tilde{\gamma}^Q(\sigma^{Q-}_{n,j})
                      \right) \\
        \st \quad& \eqref{eq:voltage_bounds}\,\text{--}\,\eqref{eq:balance_slacks_nn},
                   \eqref{eq:drop_control}, \eqref{eq:v_regulator_diff}: k = j, \\
            & 0 \leq \rho^+_{g,j},\, \rho^+_{g,j} \cdot (\overline{P}_g - p_{g,j}) \leq 
                \epsilon \cdot A_g \overline{\Delta} \cdot (\overline{P}_g - \underline{P}_g)  ~~ 
                \forall g \in \mathcal{G}^j, \\
            & 0 \leq \rho^-_{g,j},\, \rho^-_{g,j} \cdot (p_{g,j} - \underline{P}_g) \leq 
                \epsilon \cdot -A_g \underline{\Delta} \cdot (\overline{P}_g - \underline{P}_g)  ~~ 
                \forall g \in \mathcal{G}^j, \\
            & 0 \leq \nu^-_{j,k}, \, \nu^-_{j,k} \cdot (\overline{Q}_g - q_{g,j}) \leq
                \epsilon \cdot (\overline{V}_{n,j} - \underline{V}_{n, k_0}) \cdot 
                (\overline{Q}_g - \underline{Q}_g) ~~ 
                \forall g \in \mathcal{G}^j(n), n \in \mathcal{N}^j, \\
            & 0 \leq \nu^+_{j,k}, \, \nu^+_{j,k} \cdot (q_{g,j} - \underline{Q}_g) \leq
                \epsilon \cdot (\overline{V}_{n,k_0} - \underline{V}_{n, j}) \cdot 
                (\overline{Q}_g - \underline{Q}_g) ~~ 
                \forall g \in \mathcal{G}^j(n), n \in \mathcal{N}^j,
    \end{aligned} \label{eq:slave}
\end{equation}
for all $j \in \mathcal{K} \setminus \{k_0\}$.

We observe that the master subproblem~\eqref{eq:master} and the contingency
subproblems~\eqref{eq:slave} share the main OPF structure but are different in two key aspects:
(\textit{i}) subproblem \eqref{eq:slave} does not take production cost into consideration, and
(\textit{ii}) subproblem \eqref{eq:slave} is further restricted by the (approximate) drop frequency
control and voltage regulator constraints. Despite these additional constraints subproblem
\eqref{eq:slave} will be feasible as long as $(\boldsymbol{p}_{k_0}, \boldsymbol{v}_{k_0})$ is
feasible for \eqref{eq:master}
to construct a feasible solution it suffices to take
$\boldsymbol{p}_j = \boldsymbol{p}_{k_0}$ and $\boldsymbol{v}_{j} = \boldsymbol{v}_{k_0}$, and
compute values for the remaining variables accordingly. Such solution would be feasible even for
$\epsilon = 0$, satisfying constraints \eqref{eq:drop_control} -- \eqref{eq:v_regulator_down}
exactly, nevertheless this solution tends to be bad due to significant activation of power balance
slacks which result in large penalties. The next subsection examines how to recover good feasible
solutions for contingency subproblems \eqref{eq:slave} satisfying constraints
\eqref{eq:drop_control} -- \eqref{eq:v_regulator_down} exactly.

\subsection{Feasibility recovery for contingency subproblems} \label{sec:feasibility_recovery}

There are two situations under which we may require a fully feasible contingency solution:
(\textit{i}) a system operator would like to study how the system would respond to a certain
contingency or (\textit{ii}) we need a better estimate for the contingency penalties because the
relaxation \eqref{eq:slave} is lose, and further decreasing $\epsilon$ is time consuming and prone
to numerical issues. We handle these situations by, first, solving the relaxed contingency
subproblem \eqref{eq:slave} with $\epsilon > 0$ and then \emph{crushing} the approximate solution
onto the exact feasible set according to the procedure outlined below.

Let $(\boldsymbol{\tilde{p}}_j, \boldsymbol{\tilde{q}}_j)$ be an solution to problem
\eqref{eq:slave} with $\epsilon > 0$ for certain $j \in \mathcal{K} \setminus \{k_0\}$, and let
\begin{equation}
    \begin{aligned}
        \min_{\substack{p,q,v,\theta \\ \sigma, \nu, \rho}} ~~& \sum_{e \in \mathcal{E}^j}
            \left(\tilde{\gamma}^S(\sigma_{e,o,j}) + \tilde{\gamma}^S(\sigma_{e,o,j})\right) + \\
         & \sum_{n \in \mathcal{N}} 
            \left(\tilde{\gamma}^P(\sigma^{P+}_{n,j}) + \tilde{\gamma}^P(\sigma^{P-}_{n,j}) +
                  \tilde{\gamma}^Q(\sigma^{Q+}_{n,j}) + \tilde{\gamma}^Q(\sigma^{Q-}_{n,j})
                  \right) \\
        \st \quad& \eqref{eq:voltage_bounds}\,\text{--}\,\eqref{eq:balance_slacks_nn} : k = j
    \end{aligned} \label{eq:slave_restricted_canvas}
\end{equation}
be a canvas for the restricted contingency subproblem, which will be modified during the crushing
procedure. Additionally, let's define the function $\Delta_j:\Omega \rightarrow \mathbb{R}$ where
$\Omega = \big\{(x, \boldsymbol{p}_{k_0}) \in \mathbb{R} \times [\boldsymbol{\underline{P}},
\boldsymbol{\overline{P}}] \big| \sum_{g \in \mathcal{G}^j: A_g>0} (\underline{P}_g - p_{g, k_0})
< x < \sum_{g \in \mathcal{G}^j: A_g > 0} (\overline{P}_g - p_{g, k_0}) \big\}$ as
\begin{equation*}
\begin{aligned}
    \Delta_j\left(x, \boldsymbol{p}_{k_0}\right) :=~& \Bigg\{\delta_j \in \mathbb{R} ~\Big|~ 
      \exists (\boldsymbol{p}_j, \boldsymbol{\rho}^+_j, \boldsymbol{\rho}^-_j) \in 
      \mathbb{R}^{3 |\mathcal{G}^j|} ~\st~ \\
    & \eqref{eq:drop_control} \text{--} \eqref{eq:drop_control_down_saturation} : k = j, ~
      \sum_{g \in \mathcal{G}^j: A_g > 0} (p_{g,j} - p_{g, k_0}) = x \Bigg\}.
\end{aligned}
\end{equation*}
In words, $\Delta_j$ returns the value of $\delta_j$ that would cause an increase/decrease of
production $x$ in contingency $j$ with respect to a base case production set point
$\boldsymbol{p}_{k_0}$. In practice, this function can be implemented with a \emph{for} loop over the
generators in $\mathcal{G}^j$, performing a line search over $\delta_k$ ($x$ can be computed in a
straightforward fashion from $\delta_k$ and $\boldsymbol{p}_{k_0}$). 

Then, we proceed as follows. We formulate the restricted (canvas) contingency subproblem
\eqref{eq:slave_restricted_canvas}. In order to satisfy the frequency drop constraints,
\eqref{eq:drop_control} -- \eqref{eq:drop_control_down_saturation}, we modify the restricted
subproblem as indicated Alg.  \ref{alg:crushing_drop}. In the same fashion, in order to satisfy the
voltage regulation constraints, \eqref{eq:v_regulator_diff} -- \eqref{eq:v_regulator_down}, we
further modify the restricted subproblem following Alg. \ref{alg:crushing_v_reg}. Finally, we solve
the modified restricted contingency subproblem to obtain a solution satisfying all contingency $j$'s
constraints in the original SCACOPF problem. This procedure can be repeated for all contingencies to
generate a complete feasible solution to SCACOPF.

\begin{algorithm}
    \linespread{0.7}\selectfont
    \SetAlgoLined
    \KwIn{base case solution $\boldsymbol{p}_{k_0}$, 
          approximate contingency solution $\boldsymbol{\tilde{p}}_j$.}
    compute $\hat{\delta}_j := \Delta_j\big(\sum_{g \in \mathcal{G}^j} 
                                           \tilde{p}_{g, j}, \boldsymbol{p}_{k_0}\big)$\;
    compute $\hat{p}_g := \max\{\min\{p_{g,k_0} + A_g \hat{\delta}_j, \overline{P}_g\},
                                \underline{P}_g\}$ for all $g \in \mathcal{G}^j$ \;
    \eIf{$\hat{\delta}_j \geq 0$}{
        let $\check{\mathcal{G}} = \{g \in \mathcal{G}^j | \hat{p}_g = \overline{P}\}$ 
            and compute $\underline{D} = \max_{g \in \check{\mathcal{G}}} 
                                            (\overline{P}_g - p_{g,k_0})/A_g$ \;
        let $\hat{\mathcal{G}} = \{g \in \mathcal{G}^j | \hat{p}_g < \overline{P}\}$ 
            and compute $\overline{D} = \min_{g \in \hat{\mathcal{G}}} 
                                            (\overline{P}_g - p_{g,k_0})/A_g$ \;
        \With{problem \eqref{eq:slave_restricted_canvas}}{
            fix $p_{g,j} = \overline{P}_g$ for all $g \in \check{\mathcal{G}}$ \;
            set bounds $\underline{D} \leq \delta_j \leq \overline{D}$ \;
            add constraints $p_{g,j} = p_{g,k_0} + A_g \delta_j ~~ \forall g \in \hat{G}$ \;
        }
    }{
        let $\check{\mathcal{G}} = \{g \in \mathcal{G}^j | \hat{p}_g > \underline{P}\}$ 
            and compute $\underline{D} = \max_{g \in \check{\mathcal{G}}} 
                                            (\underline{P}_g - p_{g,k_0})/A_g$ \;
        let $\hat{\mathcal{G}} = \{g \in \mathcal{G}^j | \hat{p}_g = \underline{P}\}$ 
            and compute $\overline{D} = \min_{g \in \hat{\mathcal{G}}} 
                                            (\underline{P}_g - p_{g,k_0})/A_g$ \;
        \With{problem \eqref{eq:slave_restricted_canvas}}{
            fix $p_{g,j} = \overline{P}_g$ for all $g \in \hat{\mathcal{G}}$ \;
            set bounds $\underline{D} \leq \delta_j \leq \overline{D}$ \;
            add constraints $p_{g,j} = p_{g,k_0} + A_g \delta_j ~~ \forall g \in \check{G}$ \;
        }
    }
    \caption{Crushing procedure for frequency drop control constraints \eqref{eq:drop_control} --
             \eqref{eq:drop_control_down_saturation}.}
    \label{alg:crushing_drop}
\end{algorithm}

\begin{algorithm}
    \linespread{0.7}\selectfont
    \SetAlgoLined
    \KwIn{base case solution $\boldsymbol{v}_{k_0}$, 
          approximate contingency solution $\boldsymbol{\tilde{q}}_j$,
          tolerance $\epsilon^q$.}
    \For{$n \in \mathcal{N}^j$}{
        compute $\hat{\eta} = \big(\sum_{g \in \mathcal{G}^j(n)} \tilde{g}_{g,j} \big) /
                              \big(\sum_{g \in \mathcal{G}^j(n)} (\overline{Q}_g - \underline{Q}_g) \big)$ \;
        \With{problem \eqref{eq:slave_restricted_canvas}}{
            \uIf{$\hat{\eta} < \epsilon^q$}{
                fix $q_{g,j} = \underline{Q}_g$ for all $g \in \mathcal{G}^j(n)$ \;
                set bounds $v_{n,k_0} \leq v_{n,j} \leq \overline{V}_{n,j}$ \;
            }
            \uElseIf{$\epsilon^q \leq \hat{\eta} \leq 1 - \epsilon^q$}{
                fix $v_{n,j} = v_{n,k_0}$ \;
            }
            \Else{
                fix $q_{g,j} = \overline{Q}_g$ for all $g \in \mathcal{G}^j(n)$ \;
                set bounds $\overline{V}_{n,j} \leq v_{n,j} \leq v_{n,k_0}$ \;
            }
        }
    }
    \caption{Crushing procedure for voltage regulators' constraints \eqref{eq:v_regulator_diff} --
             \eqref{eq:v_regulator_down}.}
    \label{alg:crushing_v_reg}
\end{algorithm}


\subsection{Recourse functions approximation and surrogate block-incremental algorithm} \label{sec:approximate_recourse}

The key idea of our approach for handling contingencies is to approximate the recourse functions
$r_k$ in the base case subproblem \eqref{eq:master} via low-dimensional surrogates that can be,
both, handled easily by nonlinear programming solvers and updated easily after a subset of the
contingencies is evaluated. Specifically, for generator contingencies $\mathcal{K}^\mathcal{G}
\subset \mathcal{K}$ we propose the surrogate
\begin{equation*}
    \tilde{r}_k\left(p_{g(k),k_0}, q_{g(k),k_0}\right) =
        P_k \cdot \left(p_{g(k),k_0}^2 + q_{g(k),k_0}^2\right)^2
\end{equation*}
where $g(k) = \mathcal{G}^{k_0} \setminus \mathcal{G}^{k}$, that is, $g(k)$ is the failing generator
in contingency $k \in \mathcal{K}^\mathcal{G}$ and $P_k \geq 0$ is a constant to be determined. This
surrogate was inspired by three observations. First, the economics of SCACOPF (generation cost)
drive the solution towards a very limited subset of the feasible space (essentially, following merit
order) where most generators that are producing do so at their upper capacity limit. Second, within
this subset, penalties for contingencies tend to decrease whenever the quantity injected by the
failing generator is decreased (in the limit, if the generator injects zero MVA to the network, the
contingency has no effect whatsoever on the network). Third, the penalty decrease with the injection
is typically very steep and local. Trying to emulate the phenomena in the prior observations, the
proposed surrogate penalizes the fourth power of the apparent power of the failing generator,
pushing the base case towards a more risk averse solution with a rapidly increasing penalty. We
compute $P_k$ from one evaluation of $r_k(\boldsymbol{p}_{k_0}, \boldsymbol{v}_{k_0})$ as $P_k :=
r_k\big(\boldsymbol{p}_{k_0}, \boldsymbol{v}_{k_0}\big) / \big(p_{g(k),k_0}^2 +
q_{g(k),k_0}^2\big)^2$, so that $\tilde{r}_k\big(p_{g(k),k_0}, q_{g(k),k_0}\big) =
r_k\big(\boldsymbol{p}_{k_0}, \boldsymbol{v}_{k_0}\big)$ and $\tilde{r}_k\big(0, 0\big) = 0$.

We follow the same reasoning for branches (power lines and transformers) and penalize  the apparent
power entering the failing branch at the terminal where it is the highest, for all branch
contingencies $k \in \mathcal{K}^\mathcal{E}$. These surrogates produce good feasible
solutions to the SCACOPF problem~\eqref{eq:voltage_bounds} -- \eqref{eq:objective} using Alg.
\ref{alg:block_incremental_algorithm} followed by a round of feasibility recovery to obtain contingency
solutions.

\begin{algorithm}
    \linespread{0.7}\selectfont
    \SetAlgoLined
    \KwIn{maximum number of iterations $T$,
          recourse penalty tolerance $\epsilon^r$}
    set $P_k := 0$ for all $k \in \mathcal{k} \setminus \{k_0\}$ \;
    \For{$t = 1:T$}{
        solve base case subproblem \eqref{eq:master} with surrogate recourse \;
        select a block of contingencies $\mathcal{K}^t \subseteq (\mathcal{K} \setminus \{k_0\})$ \;
        solve contingency subproblem \eqref{eq:slave} for all $k \in \mathcal{K}^t$ \;
        update $P_k := r_k\big(\boldsymbol{p}_{k_0}, \boldsymbol{v}_{k_0}\big) /
                       \big(p_{g(k),k_0}^2 + q_{g(k),k_0}^2\big)^2$ 
            for all $k \in \mathcal{K}^t \cap \mathcal{K}^\mathcal{G}$ \;
        update $P_k := r_k\big(\boldsymbol{p}_{k_0}, \boldsymbol{v}_{k_0}\big) /
                       \big(p_{e(k),i,k_0}^2 + q_{e(k),i,k_0}^2\big)^2$ 
            for all $k \in \mathcal{K}^t \cap \mathcal{K}^\mathcal{E}$,
            where $i = \argmax_{j \in \{o(e), d(e)\}} p_{e(k),j,k_0}^2 + q_{e(k),j,k_0}^2$ \;
        \If{\eqref{eq:slave} evaluated \& 
            $r_k\big(\boldsymbol{p}_{k_0}, \boldsymbol{v}_{k_0}\big) < \epsilon^r$ 
            for all $k \in \mathcal{K}^t$}{break \;}
    }
    return current base case solution \;
    \caption{Surrogate block-incremental algorithm for SCACOPF.}
    \label{alg:block_incremental_algorithm}
\end{algorithm}

Algorithm \ref{alg:block_incremental_algorithm} is easily parallelizable (evaluation of contingency
subproblems), it allows for (almost) arbitrary selection of blocks of contingencies to evaluate
at each iteration, and it allows for asynchronous updates, that is, in a parallel computing setting we
can solve the base case continuously and update the surrogates as new contingency evaluations are
completed. On the selection of contingency blocks, reasonably, we would like to start evaluating
contingencies that are likely to have a large impact on the problem. How to find these contingencies
quickly and limit their impact early on is the topic of the next subsection.


\subsection{Contingency pre-screening}
\label{sec:prescreening}

Our initial approach to ranking contingencies and selecting a subset of likely impactful
contingencies was based on machine learning methodologies. This approach, whereas promising in general,
failed to provide a good ranking or classification of contingencies, its main problem being the
difficulty of generalizing machine learning models to networks that have not been seen before in the
training datasets \citep{Yang2019}. Other methodologies from the literature required extensive
calculations, which were not practical for near-real-time settings \citep{Capitanescu2007,
Fliscounakis2013}. Instead, we use intuitive selection rules, cheap to compute and with good
generalization: we select the $S^\mathcal{G}$ largest generator failures and $S^\mathcal{E}$ largest
branch failures (in terms of the elements' capacities), with $S^\mathcal{G} + S^\mathcal{E}$ smaller
or equal to the number of parallel processors available for contingency evaluations. This procedure
tends to produce some false positives (contingencies that look impactful but are not) but almost no
false negative (contingencies that look harmless but are impactful), both situations easily
remediable by exploiting parallel computing.


\input{sec_computational.tex}


\section{Numerical results} \label{sec:numerical_results}

During the course of Challenge 1, the performance of gollnlp was evaluated independently by the
organizers of the GO competition using synthetic and real SCACOPF instances. In this section we
focus on two aspects of the results of those experiments. First, we explore the qualitative
differences between an ACOPF solution and the SCACOPF solution found by our method, for a realistic
synthetic instance, and explain how such SCACOPF solution is achieved. Second, we provide a summary
of the comparative performance of gollnlp in the \emph{Final Event} of Challenge 1.

\subsection{ACOPF vs gollnlp's SCACOPF solution}

We compare two approaches for scheduling power system operations: (\textit{i}) \emph{ACOPF}, that
is, we solve the plain ACOPF problem and then evaluate its performance including contingencies, and
(\textit{ii}) \emph{SCACOPF}, \textit{e.g.}, solving SCACOPF using gollnlp. For this comparison, we use the
synthetic version of the Western Electricity Coordinating Council (WECC) developed by
\cite{wecc10k}, which we augmented with frequency drop control parameterd derived from Challenge 1
instances, stricter branch capacities based on simulations over 12 representative load cases, and
randomly generated contingencies. The resulting instance has 10\,000 buses, 12\,707 branches, 2\,485
generators, 10\,000 contingencies, and 4\,899 consumers totalling 133.1 GW of load.

\begin{figure}
\centering
\begin{minipage}{.5\textwidth}
    \centering
    {\small (a)}
\end{minipage}%
\begin{minipage}{0.5\textwidth}
    \centering
    {\small (b)}
\end{minipage}
\begin{minipage}{.5\textwidth}
    \centering
    \includegraphics[width=\textwidth]{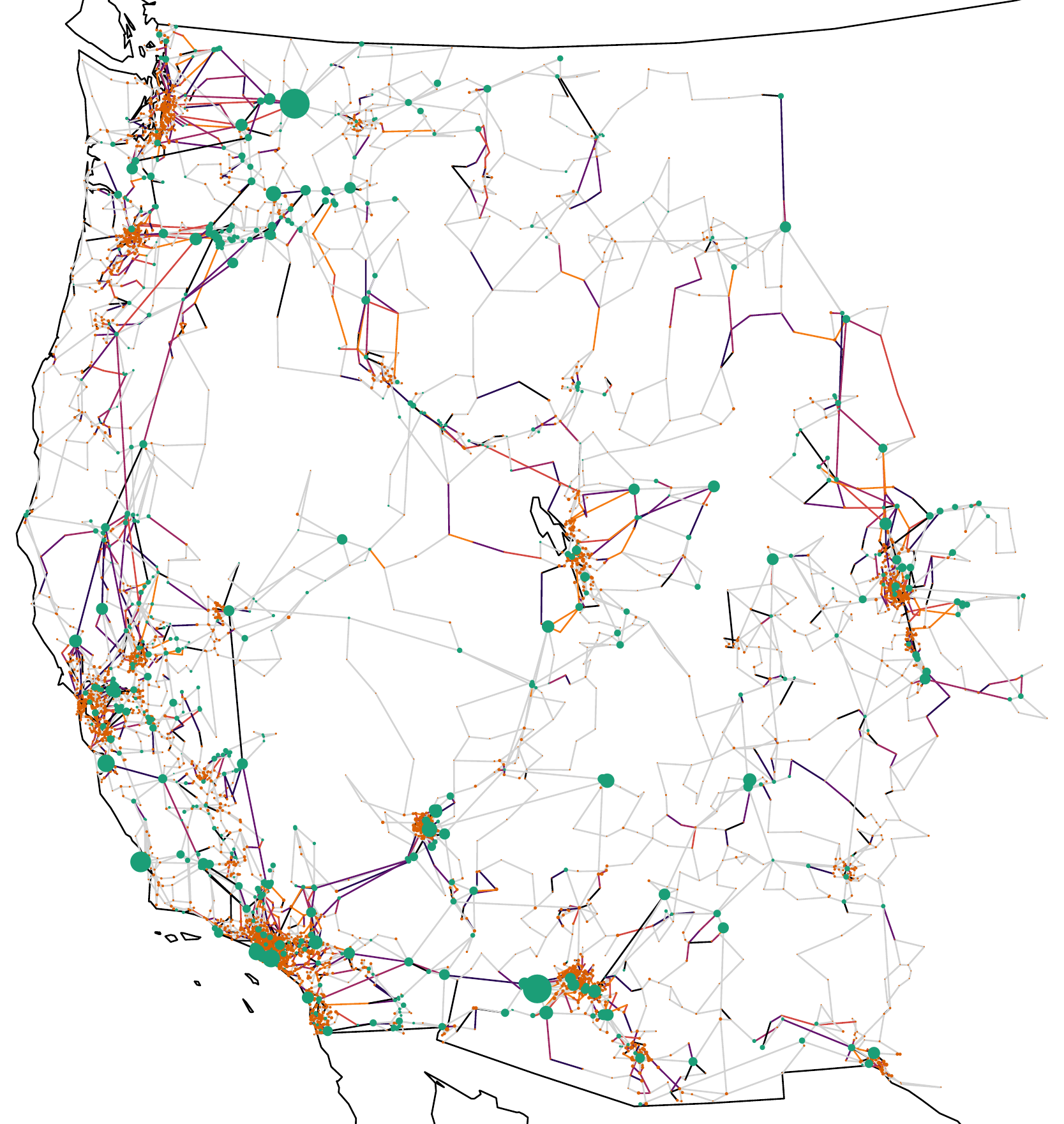}
\end{minipage}%
\begin{minipage}{0.5\textwidth}
    \centering
    \includegraphics[width=\textwidth]{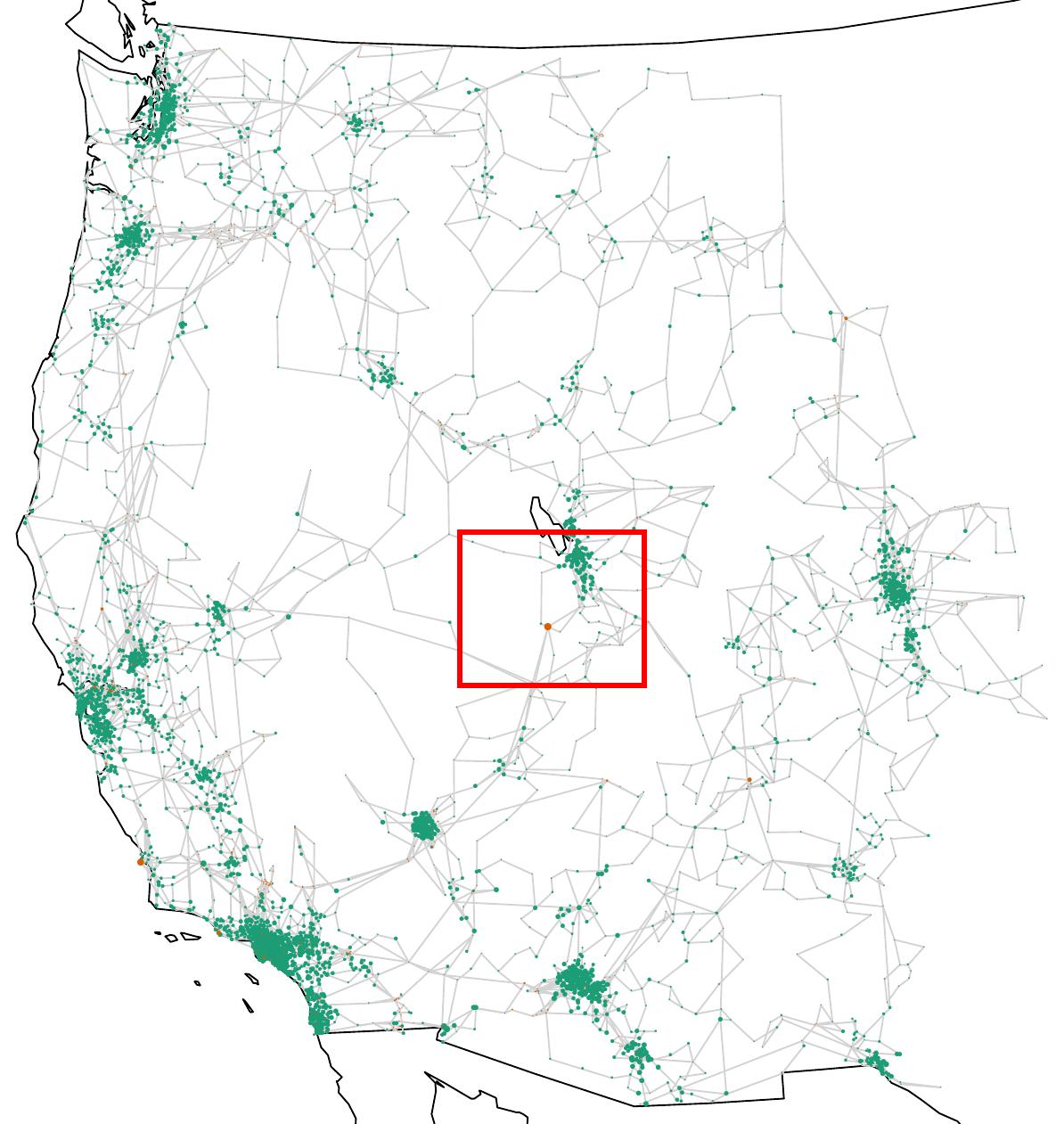}
\end{minipage}\vspace{1mm}
\caption{Synthetic WECC network color coded with \emph{ACOPF} and \emph{SCACOPF} base case
solutions. Part (a) presents the \emph{ACOPF} base case solution. Vertices correspond to
substations, colored green for substations injecting active power (generators) and orange for
substations withdrawing active power, with their area being proportional to the injected/withdrawn
active power. Edges corresponds to groups of lines between substations. They are colored according
to remaining power margin at the solution, the darker the color the smaller the margin to
congestion. Part (b) presents the injection/withdraw differences between \emph{SCACOPF} and
\emph{ACOPF}. For vertices, green indicates an increase in injection, orange indicates a decrease in
injection, and the area of each vertex is proportional to the injection difference. No information
encoded for edges. The red box delimits an area of interest, depicted in detail in Fig.
\ref{fig:zomm_box_contingency}.}
\label{fig:wecc_acopf_difference}
\end{figure}

\begin{figure}
\centering
\begin{minipage}{.33\textwidth}
    \centering
    {\small (a)}
\end{minipage}%
\begin{minipage}{0.33\textwidth}
    \centering
    {\small (b)}
\end{minipage}
\begin{minipage}{0.33\textwidth}
    \centering
    {\small (c)}
\end{minipage}\vspace{2mm}
\begin{minipage}{0.33\textwidth}
    \centering
    \includegraphics[width=.8\textwidth, frame]{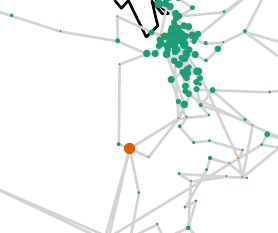}
\end{minipage}%
\begin{minipage}{0.33\textwidth}
    \centering
    \includegraphics[width=.8\textwidth, frame]{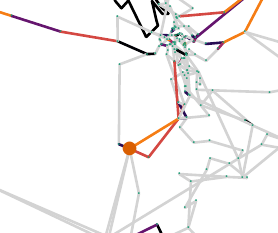}
\end{minipage}%
\begin{minipage}{0.33\textwidth}
    \centering
    \includegraphics[width=.8\textwidth, frame]{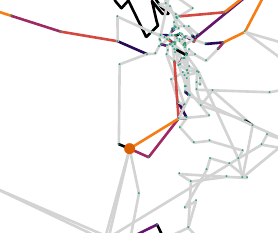}
\end{minipage}\vspace{1mm}
\caption{Differences in base case and contingency solution between \emph{ACOPF} and \emph{SCACOPF}
for synthetic WECC network. Part (a) zooms in to the area delimited with a red box in Fig.
\ref{fig:wecc_acopf_difference}.(c). Part (b) presents the post-contingency solution of \emph{ACOPF}
for the failure of a transformer at substation ``DELTA 2'' at the center of the plot. Orange vertices
indicate positive active imbalances (production shedding) green vertices indicate negative
imbalances (load shedding), with their area being proportional to the imbalance amount. Edges
colored as in Fig. \ref{fig:wecc_acopf_difference}.(a). Part (c) presents the post-contingency
solution of \emph{SCACOPF} for the same contingency, color coded in the same fashion as (b).}
\label{fig:zoom_box_contingency}
\end{figure}

Figure \ref{fig:wecc_acopf_difference} presents a depiction of the base case solution of \emph{ACOPF}
(a) and the difference on the base case solutions between \emph{SCACOPF} and \emph{ACOPF} (b).
\emph{ACOPF} schedules production minimizing base case operating cost only, thereby, setting the
cheapest generator to operate at their capacity, which tends to concentrate power generation.
\emph{SCACOPF}, in contrasts, spreads out production  -- as can be seen in Fig.
\ref{fig:wecc_acopf_difference}(b) -- in order to hedge against failures. gollnlp achieves this
hedging strategy starting from the \emph{ACOPF} solution and progressively adding approximate
recourse terms penalizing activity on elements that might fail. The mechanism at work can be
appreciated in Fig. \ref{fig:zoom_box_contingency}. The \emph{ACOPF} solution evaluated for a
transformer contingency shows a large imbalance (Fig. \ref{fig:zoom_box_contingency}(b)) because a
generator behind the transformer is dispatched at capacity in the base case and cannot reduce its
power output post-contingency due to frequency drop control constraints. gollnlp adds a recourse
term penalizing the flow over the failing transformer, which in turn rewards a reduction in the base
case dispatch of the generator behind the transformer in the \emph{SCACOPF} solution with respect to
\emph{ACOPF} solution (Fig.  \ref{fig:zoom_box_contingency}(a)).  As a result of this dispatch
reduction, the imbalance in the post-contingency state of \emph{SCACOPF} decreases significantly
(more than 100MW) with respect to that of \emph{ACOPF}, as can be seen in Fig.
\ref{fig:zoom_box_contingency}(c). This example illustrates how our non-convex relaxation of the
coupling constraints along with our recourse approximation serve to effectively characterize and
hedge against contingencies in SCACOPF.

\subsection{Benchmarking}

Our approach was benchmarked against other competing approaches in the Final Event of Challenge 1
\citep{ch1finalevent}, over 17 synthetic network models with 20 scenarios each\endnote{The Final
Event also included 3 industry network models with 4 scenarios each. We do not include statistics on
those networks in the paper as detailed results on those networks remain confidential.}, 10
scenarios corresponding to \emph{Division 1} with a time limit of 10 minutes, and the other 10
scenarios corresponding to \emph{Division 2} with a time limit of 45 minutes. All competing
approaches ran on the same platform, consisting of 6 nodes equipped with 2 Intel Xeon E5-2670 v3 (24
cores per node) and 64 GB of RAM per node, intended to represent the computational resources
currently in use by power system operators. We would like to remark that gollnlp is not bound to
these hardware specifications and it has been deployed from personal
laptops up to hundred of nodes of high-end supercomputers, such as Summit at Oak Ridge National
Laboratory.

\begin{figure}
\centering
\includegraphics[scale=.65]{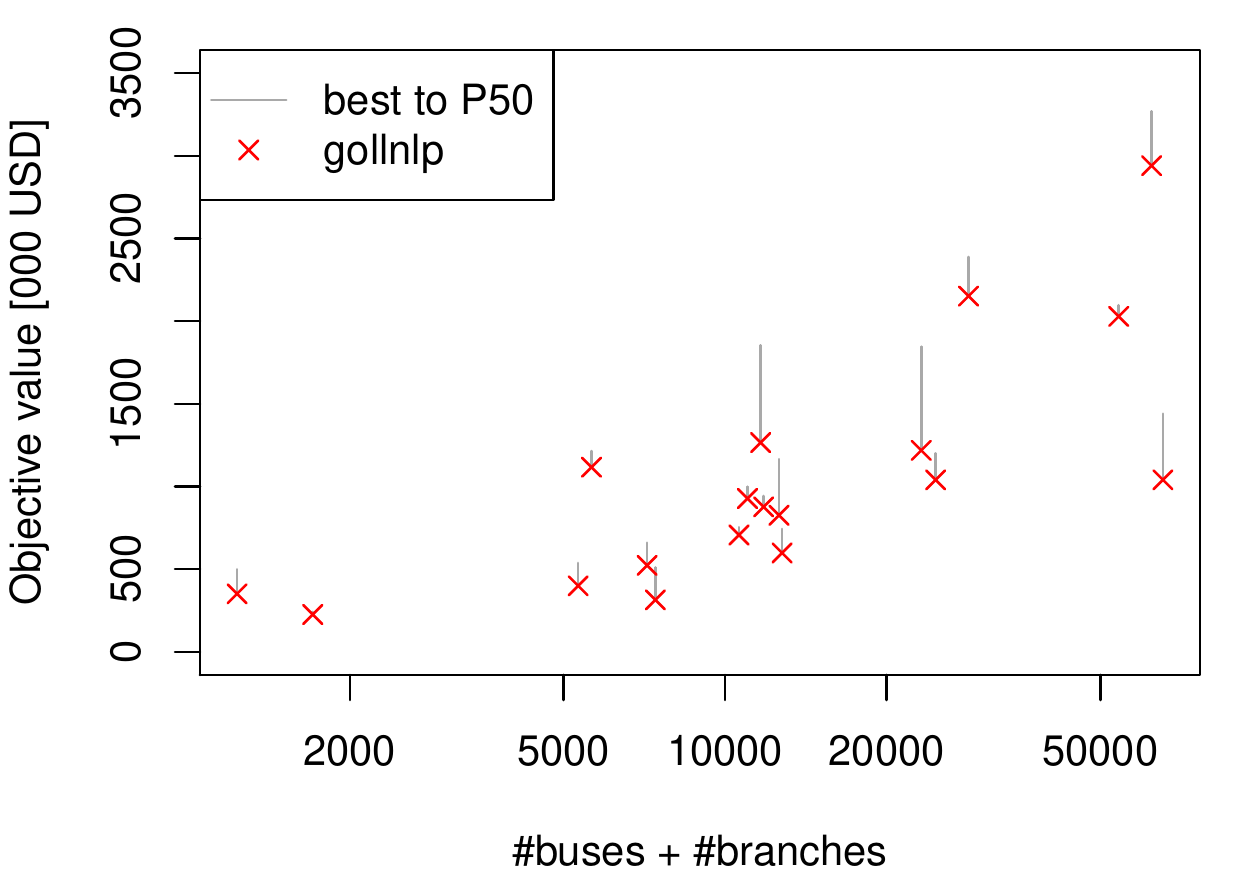}
\caption{Comparative performance of gollnlp in Division 1 of Challenge 1. Plot presents ranges of
geometric mean objective value for every synthetic network tested in the Final Event (lower is better).}
\label{fig:division1_performance}   
\end{figure}

\begin{figure}
\centering
\includegraphics[scale=.65]{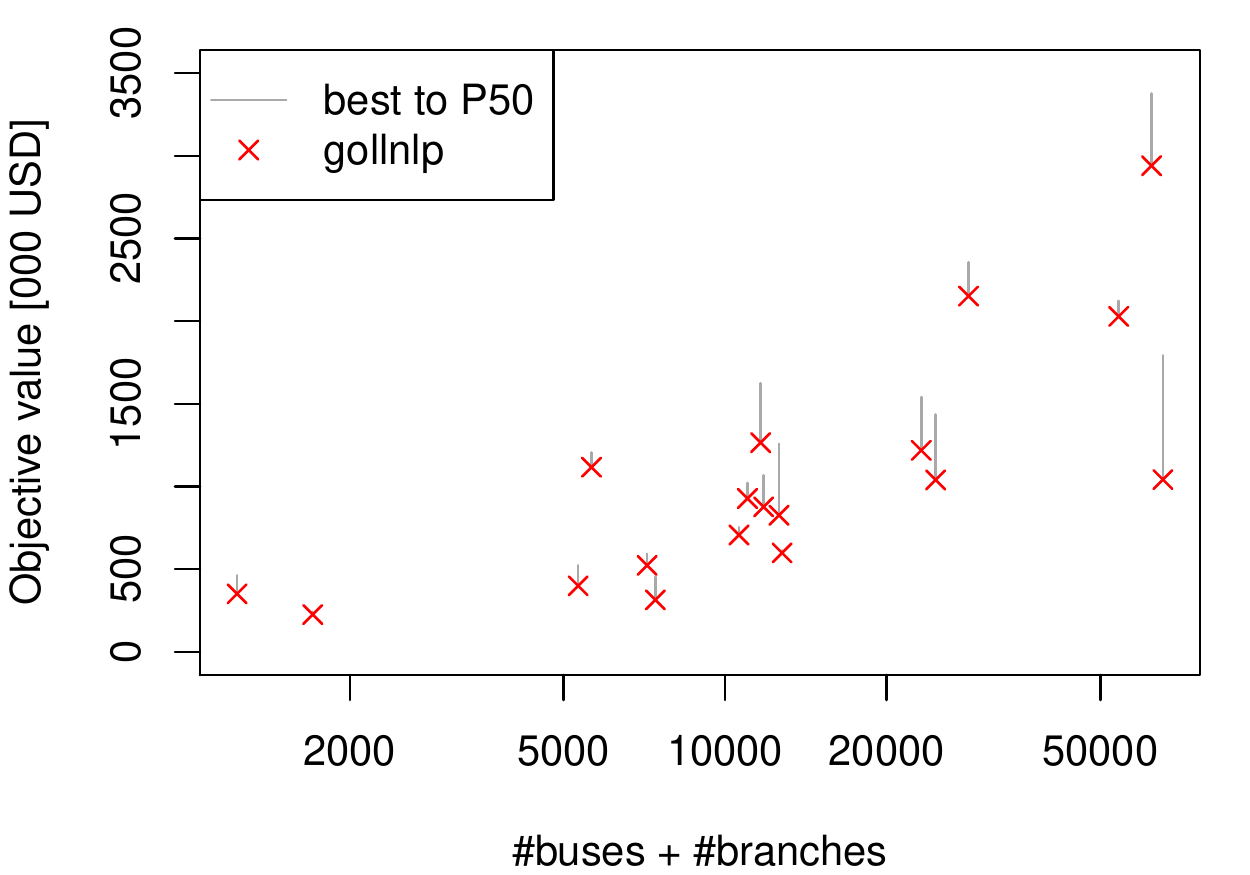}
\caption{Comparative performance of gollnlp in Division 2 of Challenge 1. Plot presents ranges of
geometric mean objective value for every synthetic network tested in the Final Event (lower is better).}
\label{fig:division2_performance}   
\end{figure}

Figures \ref{fig:division1_performance} and \ref{fig:division2_performance} present our aggregate
performance for each network in Division 1 and Division 2, respectively, against the best 50\%
objective values in each division for each network, among those able to obtain feasible solutions.
As can be observed, our approach consistently obtained the best or very close to the best known
feasible solution across a wide range of system sizes and in both Division 1 and Division 2.

\section{Conclusions} \label{sec:conclusions}

We have presented a computational approach for obtaining good feasible solutions to the SCACOPF
problem at real-world scale and under real-world time limits. Our method relies on state-of-the-art
non-linear programming algorithms and it uses non-convex relaxations of complementarity constraints,
two-stage decomposition with sparse approximations of recourse terms, and asynchronous parallelism,
in order to compute SCACOPF base case solutions that hedge against contingencies. The method was
benchmarked independently and found to consistently produce high quality solutions across a wide
range of network sizes and difficulty, thereby providing an effective complement to the extensive
literature on bounding via convex relaxations.


Future  work will investigate the efficient use of emerging computing
hardware accelerators (\textit{e.g.}, GPUs) for SCACOPF problems, as well as its global solution
using branch-and-bound techniques that integrate convex relaxations and strengthening with the
present work.

\theendnotes

\section{Acknowledgments}

This work performed under the auspices of the U.S. Department of Energy by Lawrence Livermore National Laboratory under Contract DE-AC52-07NA27344. 

\bigskip

\fbox{
\parbox{6.5in}{
\scriptsize
This document was prepared as an account of work sponsored by an agency of the United States
government. Neither the United States government nor Lawrence Livermore National Security,
LLC, nor any of their employees makes any warranty, expressed or implied, or assumes any legal
liability or responsibility for the accuracy, completeness, or usefulness of any information,
apparatus, product, or process disclosed, or represents that its use would not infringe privately
owned rights. Reference herein to any specific commercial product, process, or service by trade
name, trademark, manufacturer, or otherwise does not necessarily constitute or imply its
endorsement, recommendation, or favoring by the United States government or Lawrence
Livermore National Security, LLC. The views and opinions of authors expressed herein do not
necessarily state or reflect those of the United States government or Lawrence Livermore
National Security, LLC, and shall not be used for advertising or product endorsement purposes.
}}

\bibliographystyle{informs2014}
\bibliography{bibliography}


\end{document}

%% file: sec_computational.tex
\section{Computational and implementation aspects} \label{sec:computational_implementation}


We designed our software framework, \emph{gollnlp}, with software re-usability in mind. To this end,
we developed a lightweight model specification library that allows a modular composition of the
various algebraic ACOPF and SCACOPF models instantiated by our solution approach. These developments
and underlying high-performance computing aspects are discussed in detail in
Section~\ref{sec:sec:model}. Furthermore, the gollnlp framework has a solver-agnostic back-end that
allows switching between different NLP solvers. At the moment we support Ipopt~\citep{Watcher2006}
and HiOp~\citep{Petra_18_hiopdecomp}. We performed a couple of customizations specific to Ipopt for
the purpose of  improving robustness and reduce computational cost. These are further elaborated in
Section~\ref{sec:sec:details}.

In order to obtain a good computational parallelization efficiency of the decomposition approach
presented in Section~\ref{sec:decomposition}, we implemented a nonblocking/asynchronous parallel
evaluation engine for performing the computational analysis of the contingency subproblems and for
building their recourse evaluations presented in
Sections~\ref{sec:subproblems}--\ref{sec:approximate_recourse}. The implementation is based on the
Message Passing Interface (MPI) and is presented in Section~\ref{sec:sec:parallel}


\subsection{Model evaluation library}\label{sec:sec:model}

The specification of SCACOPF problem \eqref{eq:voltage_bounds} -- \eqref{eq:objective}  is done via
lightweight  C++ library that evaluates the objective and constraints functions, and their first and
second order derivatives, as required by the underlying NLP solvers. Both the functions and their
derivatives defining the problem \eqref{eq:voltage_bounds} -- \eqref{eq:objective} are
``hand-coded'' in C++. The specification library is modular in the sense it can assemble and
manipulate generic NLPs such as ACOPF, as well instantiate the NLPs in the low level format of
supported NLP solvers, from existing objective and constraint C++ classes. The objective can be
formed by an arbitrary number of additive objective terms, such as generation cost, penalties,
various recourse approximations, primal regularization terms, etc. The constraints of the NLP is
simply a list of (sparse) constraints of blocks of constraints. The library supports instantiating
SCACOPF problems containing a base case ACOPF block and multiple contingencies blocks and is aware
of the underlying problem structure; however, our computational methodology does not use this
feature since our optimization-based decomposition approach instantiates the first-stage master and
second-stage contingency blocks as separate problems (with the contingency subproblem being
parameterized by the first-stage's coupling variables).

Ipopt and HiOp NLP solvers require function values and gradient values as double arrays. The Jacobian of the constraints and the Hessian of the Lagragian are evaluated as sparse  matrices in triplet format. The assembly and repeated evaluation of the sparse derivatives throughout the NLP solver iterations  may easily take a high computational cost because of the high cost associated with (random) lookups (\textit{i.e.}, locating and accessing a given entry) in a sparse matrix. These lookups occur every time when a objective or constraint blocks update the derivative matrix with their (sparse) contribution and have non-trivial complexity (meaning greater than $O(1)$) in general. We circumvent this by adding a one-time preprocessing step that for each objective or constraint block builds an array of pointers to the block's nonzeros in the sparse derivative. During this preprocessing step, we sort the nonzeros based on their row and column indexes and remove duplicate entries. Consequently, the time complexity of the preprocessing step is $O(\zeta log \zeta)$ and the space requirements are $\zeta$, where $\zeta$ is the sum of the number of nonzeros of each objective or constraint block. On the other hand, we achieve $O(1)$ random access to the nonzero entries of the sparse derivative; this enables derivative evaluation at O($\zeta$) time complexity and amortizes the cost of the preprocessing step over over the course of NLP solver iterations. As a result, the evaluation of NLP consistently requires a very small fraction ($<1\%$) of the total computational cost.

\subsection{Parallel computations}\label{sec:sec:parallel}
As mentioned in Section~\ref{sec:decomposition}, the contingency subproblems~\eqref{eq:slave} can be
solved independently, in parallel, for a given solution of the master subproblem~\eqref{eq:master}.
The obvious parallelization approach, usually known as \emph{master}--\emph{worker} parallelization,
is to dedicate a small subset of the computing units to the first-stage subproblem and the rest to the
contingency subproblems. Each time all the contingencies are evaluated, the \emph{workers} send the
recourse approximations to the \emph{master} unit, which in turn resolves the first-stage problem and
dispatches the solution to the \emph{worker} units, effectively starting a new iteration of the
algorithm.

This above simple parallel programming model suffers from a couple important drawbacks.
Computational load imbalance is significant, being caused by the considerable variation across the
wall-clock (solution) times for the contingency subproblems in the \emph{worker} units. This is
especially pervasive when end-of-iteration synchronization between \emph{master} and \emph{worker}
units is enforced. An extreme, but quite frequent, consequence of this limitation manifests in a
considerable number  of \emph{worker} computing units sit idle for extended amount of time waiting
for a couple tardy \emph{workers} to finish contingency evaluations. The second limitation of a
simple \emph{master}--\emph{worker} model would be caused by the fact  that all \emph{worker} units
sit idle while the \emph{master} unit solves the first-stage subproblem. Finally, the \emph{master}
unit(s) performs significant bookkeeping, file operations, and communication with the
\emph{workers}; as a result, the \emph{master} unit may use only a subset of its compute cycles to
solve the first-stage subproblem. We remark that this results in a domino effect, as it will delay not
only the first-stage solution but will also keep all the \emph{workers} idle for longer.

We address the above limitations by using a more involved
\emph{master}--\emph{solver}--\emph{worker} parallel programming model, in which, one of the compute
units is solely assigned to solving the first-stage subproblem (hence the name \emph{solver}). This
allows the \emph{master} unit to be dedicated to algorithm control loops, file input/out, and
communication. We employ a distributed memory programming model that results in the \emph{master},
\emph{solver}, and \emph{workers} units to run in independent processes; MPI is used for
(inter-process) communication. Furthermore, the communication uses only \textit{nonblocking} MPI
calls, more specifically, only \verb~MPI_Isend~ and \verb~MPI_Irecv~. While this choice complicates
the implementation significantly because the completion of each message needs to be checked
periodically, it has the important benefit of overlapping communication with computations. This is
especially relevant to \emph{master}, which otherwise will be unable to timely handle the algorithm
control loops and the communication in the same time. 

The \emph{worker} and the \emph{solver} units engage in communication only with the \emph{master}.
The \emph{master} decides the apportionment of contingencies subproblems to \emph{worker} processes
and schedules only one contingency at a time per \emph{worker} process, which is paramount to ensure
good load balancing across the \emph{workers}. The contingency subproblems are initially prioritized
based on the rankings found in the pre-screening of Section~\ref{sec:prescreening}; later on, the
contingencies found to have high penalties are rescheduled for reevaluation/re-approximation. The
\emph{master}--\emph{worker} conversation comprises of predefined basic actions: (i) \emph{master}
asks \emph{worker} to evaluate a contingency, (ii) \emph{worker} provides contingency recourse
approximation to the \emph{master}, and (iii) initialization and finalization. On the other hand,
the \emph{master}--\emph{solver} conversation is quite simple and it involves only solve ``start''
and ``completion'' actions. In general, the \emph{master} instructs the \emph{solver} process to
reevaluate first-stage subproblem once it receives a few high-penalty  contingency recourse
approximations from the \emph{workers}. We remark that the \emph{solver} process solves the base
case problem without forcing the \emph{master} or \emph{worker} processes to sit idle; instead,
these processes can continue evaluating and approximating contingencies. Under a strict
computational budget, this approach can increase the number of evaluated contingencies by more than
$250\%$ over plain \emph{master}--\emph{worker} parallelization paradigms. The downside of our
approach is that the \emph{master} subproblem is updated with recourse approximations asynchronously
and possibly is using outdated recourse approximations built around first-stage set points from
previous iterations. The convergence of such asynchronous \emph{solver}--\emph{worker} scheme is
problematic, especially for non-convex problems such as SCACOPF. The use of recourse surrogates
based on the fourth power of the apparent power (instead of using similar quadratic surrogates) has
shown a considerable stabilization effect on the convergence the asynchronous algorithm in our
extensive numerical experiments. We believe that this is likely caused by the fact that the steep
fourth order surrogates serve as an effective proximal point term on the first-stage variables. Future
work will be dedicated to further investigate and develop globally convergent variants of this
asynchronous \emph{solver}--\emph{worker} method.

\subsection{Other implementation details}\label{sec:sec:details}

Since a great reduction of computational cost is obtained from warm-starting, we elaborate on our
experience with warm-starting interior-point methods when resolving the first-stage and the
contingency subproblems. In our computational scheme, the first-stage subproblem can be warm-started
efficiently, for $50-80\%$ reduction in the computational cost, by reusing the primal-dual solution
occurring during the previous first-stage solve for a log barrier parameter around  $10^{-5}$. For
the contingency subproblems, the best overall performance is obtained when only a primal restart is
done using the solution from the first-stage problem. For these subproblems primal-dual restarts
tends to give less benefits, regardless of using solutions from the previously solved contingency
subproblem or from the first-stage subproblem. This is mainly because a nontrivial number of
contingencies converge very slowly when using primal-dual restarts. To some extent, this also
applies when solving the feasibility recovery  subproblems, described in Section
\ref{sec:feasibility_recovery}.

Given the large size and the considerable number of contingencies of the SCACOPF problem solved,
the number of NLP subproblems that needs to be solved can be as high as 100\,000, each of
2\,000\,0000 variables and constraints. In this large-scale computational regime, the possibility of
unexpected software failures cannot be dismissed. Surprisingly, we experienced only a couple of
nontrivial issues because Ipopt is extremely robust. A repeated issue, among these, was an
occasional solver divergence or convergence to a locally infeasible point when solving the
contingency subproblems.  The issue tends to occur when solving the feasibility recovery contingency
subproblems. We could not find a sound explanation for it since the problems are provably feasible
and we alleviate these issues by adding small quadratic regularization terms in the objective for
the voltages and angles. By far the most troubling issue we ran into was that MPI worker processes
occasionally seemed to hang and -- on rare occasions -- to crash with \verb~SIGBUS~ or \verb~SIGSEGV~
signals, depending on the platform and/or compilers. This issue is aggravated because a crash in one
MPI process would cause all the MPI processes to stop execution as per the MPI standard. The
apparent source of the hang or crash was the combination of Ipopt and the HSL MA57 linear solver, but we
could not find the cause of the abnormal runtime behavior. We addressed this with a three-pronged
approach.  First we modified Ipopt to throw an exception when the linear solver unexpectedly
requires a large memory reallocation; we handled this instructing Ipopt to switching to MA27 linear
solver and resolve the problem. The second freeze and/or crash mitigation approach was to set an
alarm (\textit{a.k.a.} a timeout \verb~SIGALRM~ signal) before each call to the MA57 faulty function
and whenever the alarm signal was generated (meaning that MA57 faulty function takes an unexpectedly
long time), we modified Ipopt to throw an exception; in which case, gollnlp instructs Ipopt to use
MA27 on resolve. Interceptions of signals associated with the runtime crash were also implemented as
last resort, to prevent MPI processes from failing.

\nocite{HSL}